\documentclass[review]{elsarticle}
\usepackage{lineno,hyperref}
\usepackage{graphicx,subfig}
\usepackage{amsmath,amssymb}
\usepackage{color}
\usepackage{float}
\usepackage{tikz}
\usetikzlibrary{intersections,calc}
\usepackage{mathtools}

\modulolinenumbers[5]

\journal{Some Journal}

\begin{document}
\newcommand{\R}{\mathbb R}
\newcommand{\N}{\mathbb N}
\newcommand{\Z}{\mathbb Z}
\newcommand{\dd}{\mathrm{d}}
\newcommand{\lam}{\lambda}
\newcommand{\bfx}{\mathbf{x}}
\newcommand{\bfy}{\mathbf{y}}
\newcommand{\bfs}{\mathbf{s}}
\newcommand{\bfz}{\mathbf{z}}
\newcommand{\bfa}{\mathbf{a}}
\newcommand{\bfb}{\mathbf{b}}
\newcommand{\hK}{\widehat{K}}
\newcommand{\tK}{\widetilde{K}}
\newcommand{\cK}{\Check{K}}
\newcommand{\wA}{\widetilde{A}}
\newcommand{\hv}{\hat{v}}
\newcommand{\T}{\mathcal{T}}
\newcommand{\I}{\mathcal{I}}
\newcommand{\PP}{\mathcal{P}}
\newcommand{\Kab}{K_{ab}}
\newcommand{\bK}{\mathbf{K}}
\newcommand{\ux}{\underbar{x}}
\newcommand{\ox}{\overline{\mathrm{x}}}
\newcommand{\vvskip}{\vspace{5pt}}
\newcommand{\gmax}{\gamma_{\max}}
\newcommand{\revise}[1]{{\color{black} #1}}

\renewcommand{\figurename}{Figure~}
\renewcommand{\baselinestretch}{1.15}

\newtheorem{theorem}{Theorem}
\newtheorem{lemma}[theorem]{Lemma}
\newtheorem{definition}[theorem]{Definition}
\newtheorem{corollary}[theorem]{Corollary}
\newtheorem{condition}[theorem]{Condition}

\begin{frontmatter}

\title{A new geometric condition equivalent \\
 to the maximum angle condition for tetrahedrons}


\author[mymainaddress]{Hiroki Ishizaka}
\ead{h.ishizaka005@gmail.com}

\author[mysecondaryaddress]{Kenta Kobayashi}
\ead{kenta.k@r.hit-u.ac.jp}

\author[mymainaddress]{Ryo Suzuki}
\ead{szkryo2531g.zky@gmail.com}

\author[mymainaddress]{Takuya Tsuchiya\corref{mycorrespondingauthor}}
\cortext[mycorrespondingauthor]{Corresponding author}
\ead{tsuchiya@math.sci.ehime-u.ac.jp}

\address[mymainaddress]{Graduate School of Science and Engineering,
 Ehime University, Matsuyama, Japan}
\address[mysecondaryaddress]{Graduate School of Business Administration,
         Hitotsubashi University, Kunitachi, Japan}

\begin{abstract}
For a tetrahedron, suppose that all internal angles of faces and
all dihedral angles are less than a fixed constant  \revise{$C$} that is
smaller than $\pi$.  Then, it is said to satisfy the maximum angle
condition with the constant \revise{$C$}.  The maximum angle condition
is important in the error analysis of Lagrange interpolation on
tetrahedrons.   This condition ensures that we can obtain an error
estimation, even on certain kinds of anisotropic tetrahedrons.  In this
paper, using two quantities that represent the geometry of tetrahedrons,
we present an equivalent geometric condition to the maximum angle
condition for tetrahedrons.
\end{abstract}


\end{frontmatter}

\linenumbers

\section{Introduction}\label{intro}
Lagrange interpolation on tetrahedrons and the associated error analysis
are important subjects in numerical analysis.  They are particularly
crucial for the mathematical theory of finite element methods.
Let $T$ be an arbitrary tetrahedron and $h_T := \mathrm{diam}T$.  Let
$k$ be a positive integer, and $\PP_k(T)$ be the set of all polynomials
defined on $T$ whose degree is at most $k$.   Let
$\I_T^k v \in \PP_k(T)$ be the Lagrange \revise{interpolant} of
$v \in W^{k+1,p}(T)$, where $W^{k+1,p}(T)$ is a usual Sobolev space with
$1 \le p \le \infty$ defined on $T$ (see \cite{KobTsu20} for the
definition of Lagrange interpolation).  To obtain an error estimate such
as
\begin{align}
  \left|v - \I_T^k v\right|_{W^{m,p}(T)} \le C h_T^{k+1-m}
   |v|_{W^{k+1,p}(T)},  \quad \forall v \in W^{k+1,p}(T)
  \label{stan-est}
\end{align}
for an integer $m$ $(0 \le m \le k)$, we need to impose a certain
geometric condition on $T$.  The constant $C$ usually depends on this
geometric condition.  Two well-known geometric conditions for
tetrahedrons are listed below.
\begin{condition}[Shape-regularity condition]
Let $\rho_T$ be the \revise{diameter of the maximum ball inscribed}
 in $T$.  If
there exists a fixed constant $\sigma$ such that
\begin{align*}
   \frac{h_T}{\rho_T} \le \sigma,
\end{align*}
then $T$ is said to satisfy the \textbf{shape-regularity condition} with
$\sigma$.  This condition is also called the
\textbf{inscribed-ball condition}.
\end{condition}

The shape-regularity condition requires that $T$ is not too
``flat'' or degenerated.  If a tetrahedron satisfies the
shape-regularity condition, it is said to be \textbf{isotropic}.
If $T$ satisfies the shape-regularity condition,
then the estimate \eqref{stan-est} holds with $C = C(\sigma)$
for $k = 1$ and $3/2 < p \le \infty$ or $k \ge 2$ and
$1 \le p \le \infty$ (see  \cite[Theorem~4.4]{BreSco} and
\cite[Theorem~3.1.5]{Cia78}).  In \cite{BraKorKri08}, 
conditions that are equivalent to the shape-regularity condition are
discussed.

Besides the shape-regularity condition, the following condition is
known for the geometry of tetrahedrons.

\begin{condition}[Maximum angle condition]
For a tetrahedron $T$, suppose that there exists a constant
$\gmax \, \revise{\in [\pi/3, \pi)}$, such that all internal angles of faces
and all dihedral angles between faces are less than or equal to $\gmax$.
\end{condition}

The maximum angle condition for tetrahedrons was introduced by
K\v{r}\'{i}\v{z}ek \cite{Kri92}.  Under the maximum angle condition,
tetrahedrons may be ``flat'' or degenerated in a certain way.  Such
tetrahedrons are said to be \textbf{anisotropic}.  Estimate
\eqref{stan-est} can be proved under the maximum angle condition with
$C = C(\gmax)$ for $k = 1$ and $2 < p \le \infty$ (see \cite{Kri92} and
\cite{Duran99}).   Note that the condition $2 < p \le \infty$ cannot be
improved in the case of anisotropic tetrahedrons
\cite{KobTsu20, She94}. 

For the finite element error analysis on anisotropic meshes, readers are
refereed to  \cite{Apel, ApelDob92}.  A useful survey on the
geometric conditions associated with triangles and tetrahedrons is
given in \cite{BraHanKorKri11}. 

Recently, a new error estimation of Lagrange interpolations on
tetrahedrons was presented \cite{KobTsu20, IshKobTsu20}.
Let $h_i$ $(i = 1, \cdots, 6)$ be the edge lengths of $T$ with 
$h_1 \le h_2 \le \cdots \le h_6 = h_T$. The volume of $T$
is denoted by $|T|$.  Let $R_T$, which represents the geometry of
tetrahedrons, be defined by
\begin{align}
    R_T := \frac{h_1 h_2 h_T}{|T|}h_T.
   \label{def-RT}
\end{align}
Then, for $v \in W^{k+1,p}(T)$, we have
\begin{align}
    \left|v - \I_T^k v\right|_{W^{m,p}(T)}
   \le  C(m,k,p) \left(\frac{R_T}{h_T}\right)^m
   h_T^{k+1-m} |v|_{W^{k+1,p}(T)}
  \label{new-est}
\end{align}
with
\begin{gather}
  \begin{cases}
      2 < p \le \infty & \text{ if } k - m = 0, \\ 
      \frac{3}{2} < p \le \infty & \text{ if } k = 1, \; m = 0,\\
      1 \le p \le \infty & \text{ if } k \ge 2 \text{ and } \; k-m \ge 1.
   \end{cases}
   \label{p-cond}
\end{gather}
Note that the constant $C(m,k,p)$ depends only on
$m$, $k$, $p$ (and the space dimension $3$).  Therefore, we can 
apply estimation \eqref{new-est} to arbitrary tetrahedrons.

\vspace{3mm}
\noindent
\textit{Remark.} In \cite{KobTsu20}, $R_T$ was defined as 
the \textit{projected circumradius} of $T$ and \eqref{new-est} was
proved.  In \cite{IshKobTsu20}, $R_T$ is redefined as \eqref{def-RT},
making the proof of \eqref{new-est} much simpler.
\footnote{See also [Erratum Corollary 2D] in 
\texttt{https://arxiv.org/abs/2002.09721}.}
It is conjectured
that $R_T$ defined by \eqref{def-RT} and the projected circumradius
of $T$ are equivalent.

\vspace{3mm}
For Lagrange interpolation on triangles, a similar estimation
to \eqref{new-est} holds by setting $R_T$ as the circumradius of $T$
\cite{KobTsu15}.  From the law of sines, we realize that a triangle $T$
satisfies the maximum angle condition (which is defined in a similar
manner) if and only if there exists a fixed constant $D$ such that
\begin{align*}
    \frac{R_T}{h_T}\, \revise{= \frac{1}{2 \sin\theta_T}}\le D,
\end{align*}
\revise{where $\theta_T$ is the maximum internal angle of $T$.
If the above inequality holds for a triangle $T$ with a fixed constant
$D$, then the triangle is said to satisfy the \textit{semiregularity}
condition with $D$ \cite{Kri91}. }

The aim of this paper is to prove the following theorem, which claims
that a similar situation holds for tetrahedrons.

\begin{theorem}[Main theorem] \label{main-thm}
Let $T \subset \R^3$ be an arbitrary tetrahedron and $R_T$ be defined by
\eqref{def-RT}.  Then, $T$ satisfies the maximum angle condition with
$\gmax \, \revise{\in [\pi/3, \pi)}$, if and only if there exists a fixed
constant $D = D(\gmax)$ such that
\begin{align}
    \frac{R_T}{h_T} \le D.
    \label{equiv-cond}
\end{align}
\end{theorem}

\begin{corollary}\label{cor1}
Let $k$ be a positive integer. Suppose that a
tetrahedron $T$ satisfies the maximum angle condition with
$\gmax \, \revise{\in [\pi/3, \pi)}$.
 Then, for the Lagrange interpolation $\I_T^k$
on $T$, the following estimate holds:
\begin{align*}
  |v - \I_T^k v|_{W^{m,p}(T)} \le C
   h_T^{k+1-m} |v|_{W^{k+1,p}(T)},
  \quad \forall v \in W^{k+1,p}(T),
\end{align*}
where $m$, $p$ are taken as in \eqref{p-cond} and the constant
$C$ depends on $m$, $k$, $p$, and $\gmax$.
\end{corollary}

\revise{
The above mentioned estimations 
(\eqref{stan-est}, \eqref{new-est}, and Corollary~\ref{cor1}) give 
upper bounds of the interpolation errors on a single tetrahedron $T$.  Error
estimations of the global Lagrange interpolation $\I_{\T_h}^k$ defined on a
simplicial mesh $\T_h$ of a bounded polyhedral domain $\Omega$ can be
obtained as
\begin{align*}
   |v - \I_{\T_h}^k v|_{m,p,\Omega} =
   \left(\sum_{T \in \T_h} |v - \I_T^k v|_{m,p,T}^p\right)^{1/p},
   \quad \forall v \in W^{k+1,p}(\Omega).
\end{align*}
Hence,  from \eqref{new-est}, we have
\begin{align*}
   |v - \I_{\T_h}^k v|_{m,p,\Omega} & \le
    C_1 \max_{T \in \T_h} \left(R_T^m h_T^{k+1-2m}\right)
     |v|_{k+1,p,\Omega} \\
   & \le
    C_1 R^m h^{k+1-2m}
     |v|_{k+1,p,\Omega},  \quad \forall v \in W^{k+1,p}(\Omega),
\end{align*}
where $h := \max_{T \in \T_h} h_T$,  $R := \max_{T \in \T_h} R_T$, and
 $C_1 = C_1(m,k,p)$ is independent of the geometry of the
tetrahedrons in $\T_h$.  If all $T \in \T_h$ satisfy the maximum angle
condition with a fixed constant $\gmax$, we have, from
Corollary~\ref{cor1},
\begin{align*}
   |v - \I_{\T_h}^k v|_{m,p,\Omega} \le
    C_2 h^{k+1-m} |v|_{k+1,p,\Omega},  \quad \forall v \in W^{k+1,p}(\Omega),
\end{align*}
where $C_2 = C_2(m,k,p,\gmax)$.
}

\section{Preliminaries}\label{prelim}

\subsection{Notation}
Let $T$ be a tetrahedron in $\R^3$ with vertices $P_1$,
$P_{2}$, $P_{3}$, and $P_{4}$.   The edge connecting $P_i$ and $P_j$ and
its length are denoted by $\overline{P_iP_j}$ and $|\overline{P_iP_j}|$,
respectively ($i,j = 1, 2, 3, 4$, $i \neq j$).

We introduce the following notation convention on $T$.  Let $F_i$ be the
face of $T$ opposite to $P_i$.  We denote  the dihedral angle between
the faces $F_i$ and $F_j$ by $\psi^{i,j}$. Note that
$\psi^{i,j} = \psi^{j,i}$. Furthermore, we denote the internal angle at
$P_j$ on $F_i$  by $\theta_j^i$, and the angle between $F_i$ and
$\overline{P_iP_j}$ by $\phi_j^i$. 

\begin{center}
Table~1. Notation convention on $T$ 
$(i,j = 1, 2, 3, 4, \; i \neq j)$. \\[6pt]
\begin{tabular}{|c|l|}\hline
  $P_i$  & the vertices of $T$. \\ \hline
  $F_i$  & the face opposite to $P_i$. \\ \hline
  $\psi^{i,j}$ 
   &  the dihedral angle between $F_i$ and $F_j$. \\ \hline
  $\theta_j^i$ 
   &  the internal angle of $F_i$ at $P_j$. \\ \hline
  $\phi_j^i$ 
   & the angle between $F_i$ and $\overline{P_i P_j}$. \\ \hline
 \end{tabular}
\end{center}

\begin{figure}[htbp]
\begin{tikzpicture}[line width = 1pt,scale=1.0]
   \coordinate [label=below:{$P_n$}](A) at (0.5,0.0);
   \coordinate [label=right:{$P_k$}](B) at (4.5,1.5);
   \coordinate [label=above:{$P_j$}](D) at (2.0,5.0);
   \coordinate [label=left:{$P_m$}](C) at (-1.0,2.0);
   \draw[line width=1pt] (A) -- (B) -- (D) -- (C) -- (A) -- (D);
   \draw[dotted] (B) -- (C);
   \draw[thin] ([shift={(A)}] 20:0.7) 
    arc[radius=0.7, start angle=20, end angle=75]
    node[right]{$\quad \theta_n^m$};
   \coordinate (G) at ($(A)!0.13!(D)$);
   \coordinate (H) at ($(A)!0.2!(C)$);
   \draw[line width=0.3pt] (G) to[bend right]
    node[above]{{$\theta_n^k$}} (H);
\end{tikzpicture}
\qquad
\begin{tikzpicture}[line width = 1pt,scale=1.0]
   \coordinate [label=below:{$P_n$}](A) at (0.5,0.0);
   \coordinate [label=right:{$P_k$}](B) at (4.5,1.5);
   \coordinate [label=above:{$P_j$}](D) at (2.0,5.0);
   \coordinate [label=left:{$P_m$}](C) at (-1.0,2.0);
   \coordinate [label=below:{\revise{$A$}}](E) at (1.5,1.3);
   \coordinate [label=below:{\revise{$B$}}](F) at ($(A)!0.55!(B)$);
   \draw[line width=1pt] (A) -- (B) -- (D) -- (C) -- (A) -- (D);
   \draw[dotted] (B) -- (C);
   \draw[line width=0.3pt] (A) -- (E) -- (D);
   \draw[line width=0.3pt] (E) -- (F) -- (D);
   \coordinate (G) at ($(A)!0.15!(D)$);
   \coordinate (H) at ($(A)!0.3!(E)$);
   \draw[line width=0.3pt] (G) to[bend left]
    node[right]{\!\raisebox{35pt}{$\phi_n^j$}} (H);
   \coordinate (P) at ($(F)!0.2!(E)$);
   \coordinate (Q) at ($(F)!0.1!(D)$);
   \draw[line width=0.3pt] (P) to[bend left]
    node[left]{\raisebox{20pt}{$\psi^{j,m}$}\!\!} (Q);
\end{tikzpicture}
\caption{Definitions of the angles on $T$.} \label{fig1}
\end{figure}
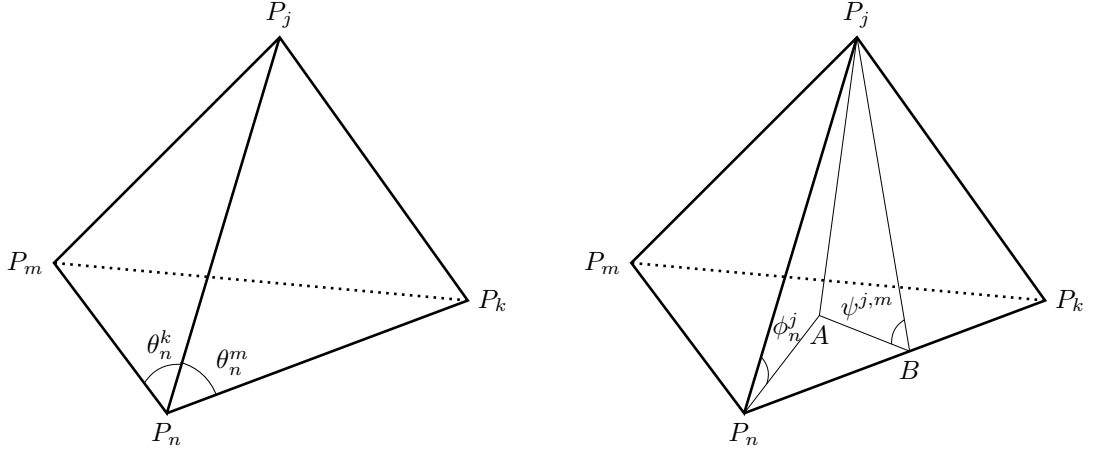

\revise{
Let $A$ and $B$ be the feet of perpendicular lines from $P_j$ to $F_j$
and from $P_j$ to $\overline{P_nP_k}$, respectively (see Figure~\ref{fig1}).
Then, we have
\begin{align*}
  |\overline{P_jP_n}| \sin \phi_n^j = |\overline{P_jA}|
   = |\overline{P_jB}| \sin \psi^{j,m}
   = |\overline{P_jP_n}| \sin \theta_n^m \sin \psi^{j,m}.
\end{align*}
A similar equation holds for $\phi_n^j$, $\theta_n^k$, and $\psi^{k,j}$.
Therefore,}
\begin{align}
  \begin{aligned}
  &  \sin \phi_n^j = \sin \theta_n^k \sin \psi^{k,j} =
    \sin \theta_n^m \sin \psi^{m,j}\revise{,} \\
  & \qquad j = 1, 2, 3, 4,  \quad m, n, k
    \in \{1,2,3,4\} \backslash \{j\}.
   \end{aligned}
  \label{twosin}
\end{align}

\subsection{Classification of tetrahedrons into two types}
As noted in \cite{Apel, IshKobTsu20, KobTsu20}, 
to deal with arbitrary tetrahedrons (including anisotropic ones)
uniformly, we need to classify tetrahedrons into two types.   Let $e_2$
be the shortest edge of $T$ and $e_1$ be the longest edge 
connected to $e_2$.  We assume that that $P_1$ and $P_2$ are the
endpoints of $e_1$, and that $e_2$ is an edge of the face
$F_4 = \triangle P_1P_2P_3$.

Consider the plane that is perpendicular to $e_1$ and intersects
$e_1$ at its midpoint.  Then, $\R^3$ is divided by this plane into two
half-spaces.  In this situation, we have two cases,
and tetrahedrons are classified as either Type~1 or Type~2 accordingly:
\begin{itemize}
 \item \textbf{Case~1}. If vertices $P_3$ and $P_4$ belong to the
same half-space, \\ \hspace{13.5mm} then $T$ is classified as Type~1.
 \item \textbf{Case~2}. If vertices $P_3$ and $P_4$ belong to
different half-spaces, \\ \hspace{13.5mm} then $T$ is classified as Type~2.
\end{itemize}
If $P_3$ or $P_4$ is on the plane, we suppose that $P_3$ and $P_4$
belong to the same half-space.   

We now introduce the following assignment of the vertices for each case.
\begin{itemize}
 \item If $T$ is Type~1, the endpoints of $e_2$ are
$P_1$ and $P_3$,  and $\alpha_2 := |\overline{P_1P_3|}$.
 \item If $T$ is Type~2, the endpoints of $e_2$ are
$P_2$ and $P_3$,  and $\alpha_2 := |\overline{P_2P_3|}$.
\end{itemize}
Define $\alpha_1 := |\overline{P_1P_2|}$ and 
$\alpha_3 := |\overline{P_1P_4|}$ for both cases.  Note that we have
implicitly assumed that $P_1$ and $P_4$ belong to the same half-space
for both cases.

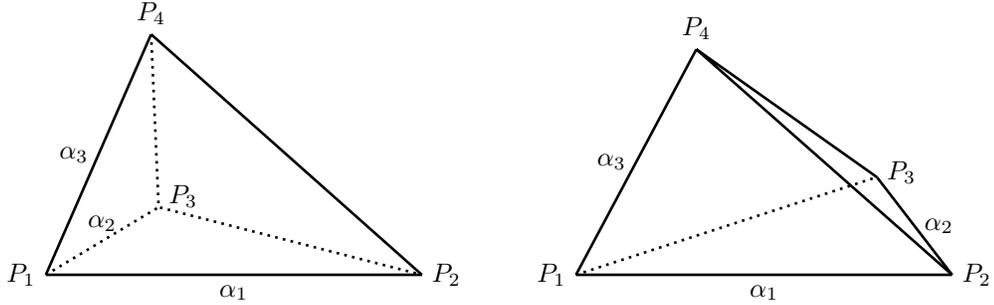
\begin{figure}[htbp]
\begin{tikzpicture}[line width = 1pt]
   \coordinate [label=left:{$P_1$}](A) at (0.0,0.0);
   \coordinate [label=right:{$P_2$}](B) at (5.0,0.0);
   \coordinate [label=above:{$P_4$}](D) at (1.4,3.2);
   \coordinate [label=right:{\raisebox{10pt}{$P_3$}}](C) at (1.5,0.9);
   \draw (A) to node[below]{$\alpha_1$}(B);
   \draw (B) -- (D);
   \draw (A) to node[left]{$\alpha_3$} (D);
   \draw[dotted] (D) -- (C) -- (B);
   \draw[dotted] (A) to node[above]{$\alpha_2$} (C);
\end{tikzpicture}
\qquad
\begin{tikzpicture}[line width = 1pt]
   \coordinate [label=left:{$P_1$}](A) at (0.0,0.0);
   \coordinate [label=right:{$P_2$}](B) at (5.0,0.0);
   \coordinate [label=above:{$P_4$}](D) at (1.6,3.0);
   \coordinate [label=right:{\raisebox{5pt}{$P_3$}}](C) at (4.0,1.3);
   \draw (A) to node[below]{$\alpha_1$}(B);
   \draw (B) to node[right]{$\alpha_2$}(C);
   \draw (C) -- (D) -- (B);
   \draw (A) to node[left]{$\alpha_3$} (D);
   \draw[dotted] (A) -- (C);
\end{tikzpicture}
\caption{Tetrahedrons of Type~1 (left) and Type~2 (right).} \label{fig2}
\end{figure}

\subsection{Another quantity that represents the geometry of $T$}
\label{HT}
For a tetrahedron $T$, we define $H_T$, which also represents
the geometry of $T$, by
\begin{align}
   H_T := \frac{\alpha_1\alpha_2\alpha_3}{|T|}h_T.
\end{align}
\revise{The values} $H_T$ and $R_T$, defined by \eqref{def-RT}, have the following
equivalence \cite[Lemma~3]{IshKobTsu20}.
\begin{lemma}
 For an arbitrary tetrahedron $T$, $H_T$ and $R_T$
are equivalent:
\begin{align*}
  \frac{1}{2}H_T \le R_T \le 2 H_T.
\end{align*}
\end{lemma}
Therefore, to prove Theorem~\ref{main-thm}, we may use $H_T$ instead of
$R_T$.

\section{Lemmas}
In this section, we prepare some useful lemmas.
In the following, we abbreviate ``maximum angle condition''
as MAC.
 
\begin{lemma}[Cosine rules on tetrahedrons]
Let $T \subset \mathbb{R}^3$ be a tetrahedron. 
Let $j = 1, 2, 3, 4$ and $\{k, m, n\} = \{1,2,3,4\} \backslash \{j\}$.
Then, we have
\begin{align}
\cos \theta^{k}_j & = \cos \theta^{m}_j \cos \theta^{n}_j
      + \sin \theta^{m}_j \sin \theta^{n}_j \cos \psi^{m,n},
   \notag \\
\cos \psi^{{n},{m}} & = \sin \psi^{m,k} \sin \psi^{n,k}
   \cos \theta^{k}_j - \cos \psi^{m,k} \cos \psi^{n,k}.
   \label{semi17d}
\end{align}
\end{lemma}
\noindent
\textit{Proof.}
See \cite{GelCot,Tod1886}.  $\square$

\begin{lemma} \label{lem53}
Let $T \subset \mathbb{R}^2$ be a triangle and let $\theta_i$
$(i=1,2,3)$  be the internal angles of $T$ with
$\theta_1 \leq \theta_2 \leq \theta_3$. If there exists
$\gmax \, \revise{\in [\pi/3, \pi)}$ such that
$\theta_{3} \leq \gmax$, then we have
\begin{align}
\sin \theta_{2}, \  \sin \theta_{3} \geq \min \left\{
   \sin \frac{\pi - \gmax}{2}, \sin \gmax \right\}.
   \label{Lem7}
\end{align}
\end{lemma}
\noindent
\revise{
\textit{Proof.} Because $\theta_1 + \theta_2 + \theta_3 = \pi$,
the assumptions yield 
\begin{align*}
 2\theta_2 \ge \theta_1 + \theta_2
 = \pi - \theta_3 \ge \pi - \gmax \quad
  \text{ and } \quad
   \frac{\pi - \gmax}{2} \le \theta_2 \le \theta_3 \le \gmax,
\end{align*}
which implies \eqref{Lem7}.   $\square$
}


\begin{lemma} \label{lem55}
For $\gamma \in \, \revise{[{\pi}/{3}, \pi)}$,
we have
\begin{align*}
\displaystyle
0 < \frac{\cos \gamma + 1}{\sin \frac{\gamma}{2} + 1} \leq 1.
\end{align*}
\end{lemma}
\noindent
\revise{
\textit{Proof.} This lemma can be proved immediately from
\begin{align*}
   \frac{\cos\gamma + 1}{\sin \frac{\gamma}{2} + 1}
  = 2 \left(1 - \sin\frac{\gamma}{2}\right), \qquad
   \frac{\pi}{6} \le \frac{\gamma}{2} < \frac{\pi}{2}, \qquad
   \frac{1}{2} \le \sin\frac{\gamma}{2} < 1. 
  \qquad \square
\end{align*}
}
\begin{lemma} \label{lem56}
Let $T \subset \R^3$ be a tetrahedron.  Suppose that $T$ satisfies the
MAC with $\gmax \, \revise{\in [\pi/3, \pi)}$.
Additionally, assume that $\theta_n^j$ is not the minimum
angle of face $F_j = \triangle P_m P_nP_k$, and
$\theta_n^j < {\pi}/{2}$, where $j = 1,2,3,4$ and
$\{m,n,k\} = \{1,2,3,4\} \backslash \{j\}$.
Then, setting $\delta$ to
\begin{align*}
\displaystyle
\sin \delta = \left( \frac{\cos \gmax + 1}
   {\sin \frac{\gmax}{2} + 1} \right)^{1/2}, \qquad
   0 <  \delta \leq \frac{\pi}{2},
\end{align*}
we have either
\begin{align}
\displaystyle
\psi^{m,j} \geq \delta, \quad \text{or} \quad
   \psi^{k,j} \geq \delta. \label{semi18}
\end{align}
\end{lemma}
\noindent
\textit{Proof.} From Lemma~\ref{lem55}, we have
\begin{align*}
\displaystyle
0 < \frac{\cos \gmax + 1}{\sin \frac{\gmax}{2} + 1} \leq 1,
\end{align*}
and we confirm that $\delta$ is well-defined.

The proof is by contradiction.  Suppose that 
\begin{align*}
0 < \psi^{m,j} < \delta \quad \text{ and } \quad 0 < \psi^{k,j} < \delta.
\end{align*}
Then, we have
$0 < \sin  \psi^{m,j} \sin \psi^{k,j} < \sin^2 \delta$ and
$1 > \cos  \psi^{m,j} \cos \psi^{k,j} > \cos^2 \delta$. 
From Lemma \ref{lem53} and the assumption, we have
\begin{align*}
\frac{\pi - \gmax}{2} \leq \theta_n^j < \frac{\pi}{2},  \qquad
 0 < \cos  \theta_n^j \le 
  \cos \left( \frac{\pi - \gmax}{2} \right)
   = \sin \frac{\gmax}{2}.
\end{align*}
Thus, we obtain
\begin{align*}
  \sin  \psi^{m,j} \sin \psi^{k,j} \cos \theta_n^j
   < \sin^2 \delta \sin \frac{\gmax}{2}.
\end{align*}
The cosine rule \eqref{semi17d} and the above inequalities yield
\begin{align*}
\cos \psi^{m,k} & = \sin \psi^{m,j}
  \sin \psi^{k,j} \cos \theta_n^{j}
   - \cos \psi^{m,j} \cos \psi^{k,j} \\
 & < \sin^2 \delta \sin \frac{\gmax}{2} - (1 - \sin^2 \delta ) \\
&= \frac{\cos \gmax + 1}{\sin \frac{\gmax}{2} + 1}
 \left( \sin \frac{\gmax}{2} + 1 \right) - 1 =  \cos
 \gmax,
\end{align*}
which contradicts the MAC: $\psi^{m,k} \leq \gmax$.  \qed

\begin{corollary} \label{coro57}
Under the assumptions of Lemma \ref{lem56}, we have
\begin{align*}
\sin \psi^{m,j} \ge C_0, \quad \text{or}
 \quad \sin \psi^{k,j} \ge C_0, \qquad
  C_0 := \min \{ \sin\delta , \sin\gmax \}.
\end{align*}
\end{corollary}

\begin{lemma} \label{lem59}
For $j = 1,2,3, 4$, let $\{m,n,k\} = \{1,2,3,4\} \backslash \{j\}$.
Let $p \in \{m, n, k\}$, and $\{q, r\} = \{m, n, k\} \backslash \{p\}$.
Suppose that there exists a positive constant $M$ 
with $0 < M < 1$ such that $\sin \phi_p^j \sin\theta_n^j\ge M$.
Then, setting $\gamma (M) := \pi - \sin^{-1} M$
$(\frac{\pi}{2} < \gamma(M) < \pi)$, the MAC with
$\gamma(M)$ is satisfied on faces $F_j$, $F_q$, $F_r$, and
$\psi^{j,q}$, 
$\psi^{j,r} \, \revise{\le} \, \gamma(M)$.
\end{lemma}
\noindent
\textit{Proof.} From the assumption, we have
\begin{align*}
 M \, \revise{\le} \,   
  \sin \phi_p^j \sin\theta_n^j \le \sin \theta_n^j
   \quad \text{ and } \quad
  M \, \revise{\le} \, \sin \phi_p^j.
\end{align*}
Hence, the definition of $\gamma(M)$ yields 
\revise{$\pi - \gamma(M) \le \theta_n^j \le \gamma(M)$.}
\revise{Because $\theta_n^j + \theta_m^j + \theta_k^j = \pi$,
we see that $\theta_m^j$, $\theta_k^j < \theta_m^j + \theta_k^j
  \le \gamma(M)$.
}
%
%
That is, the MAC with $\gamma(M)$ is satisfied on
face $F_j = \triangle P_mP_nP_k$.

Moreover, it follows from
\eqref{twosin} that 
\begin{align*}
  M & \, \revise{\le} \, \sin \phi_p^j = \sin \theta_p^q \sin \psi^{q,j} 
      = \sin \theta_p^r \sin \psi^{r,j} \\
     & \le \sin \theta_p^q, \; \sin \theta_p^r, \;
        \sin \psi^{r,j}, \; \sin \psi^{q,j} 
\end{align*}
By the same reasoning, we find that the MAC
with $\gamma(M)$ is satisfied on faces $F_q$ and $F_r$, and
$\psi^{j,q}$, $\psi^{j,r} \,\revise{\le} \, \gamma(M)$.  $\square$

\section{Proof of Theorem~\ref{main-thm}}
In this section, we prove Theorem~\ref{main-thm}.
As explained in Section~\ref{HT}, we may and will use
$H_T$ instead of $R_T$ in the proof.
We divide the proof into four cases.

\subsection{Type~1: Proof of ``MAC implies \eqref{equiv-cond}''}
\label{Case1}
First, we suppose that $T$ is of Type~1 and satisfies the MAC
with $\gmax$, $\pi/3 \le \gmax < \pi$.
Because $|T| = \frac{1}{6} \alpha_1 \alpha_2 \alpha_3
\sin \theta_1^4  \sin \phi_1^4$, we have
 \begin{align*}
   \frac{H_T}{h_T} = \frac{\alpha_1 \alpha_2 \alpha_3}{|T|}
     = \frac{6}{\sin \theta_1^4  \sin \phi_1^4}.
\end{align*}
From the definition of Type~1, we realize that
$\theta_2^4 \le \theta_1^4 \le \theta_3^4$, that is,
$\theta_3^4$ and $\theta_2^4$ are the maximum and minimum angles
of face $F_4 = \triangle P_1 P_2 P_3$, respectively. 
Thus, it follows from Lemma \ref{lem53} that
\begin{align*}
 \frac{\pi - \gmax}{2} \leq \theta_1^4 \leq \gmax, \quad  \sin
 \theta_1^4 \geq \min \left \{ \sin \frac{\pi - \gmax}{2}, \sin
 \gmax \right \} =: C_1.
\end{align*}
Additionally, we may apply Lemma~\ref{lem56} to $\theta_1^4$ and $F_4$,
and find that either $\psi^{2,4} \ge \delta$ or
$\psi^{3,4} \ge \delta$, where
$\delta = \delta(\gmax)$, $0 < \delta \le \pi/2$ is defined as
\begin{align}
  \sin \delta = \left( \frac{\cos \gmax + 1}
   {\sin \frac{\gmax}{2} + 1} \right)^{1/2}.
   \label{def-delta}
\end{align}

Suppose that $\psi^{2,4} \ge \delta$. By Corollary~\ref{coro57}
and \eqref{twosin}, we have 
\begin{align*}
  \sin \phi_1^4 = \sin \theta^2_1 \sin \psi^{2,4}
  \ge C_0 \sin \theta^2_1,
\end{align*}
where $C_0$ is the constant defined in Corollary~\ref{coro57}.
By the definition of Type 1, $\theta_1^2$ is not
the minimum angle of $F_2 = \triangle P_1 P_3 P_4$, and therefore, we
have
\begin{align*}
  \frac{\pi - \gmax}{2} \le \theta_1^2 \leq \gmax,
  \quad  \sin \theta_1^2 \ge  C_1.
\end{align*}
Thus, we obtain $\sin \phi_1^4 \geq C_0 C_1$.

Next, suppose that $\psi^{3,4} \geq \delta$.
Replacing $\psi^{2,4}$, $\theta_1^2$, and $F_2$ with
$\psi^{3,4}$, $\theta_1^3$, and $F_3$ in the above argument,
we obtain $\sin \phi_1^4 \geq C_0 C_1$ in the same manner.

Gathering the above results, we conclude that
\begin{align*}
  \frac{H_T}{h_T} = \frac{6}{\sin \theta_1^4  \sin \phi_1^4}
  \le \frac{6}{C_0 C_1^2} =: D
\end{align*}
in both cases, that is, \eqref{equiv-cond} holds.

\subsection{Type~1: Proof of ``\eqref{equiv-cond} implies MAC''}
\label{Case2}
Now, we suppose that $T$ is of Type~1 and 
\begin{align*}
  \frac{H_T}{h_T} = \frac{\alpha_1 \alpha_2 \alpha_3}{|T|}
   = \frac{6}{\sin \theta_1^4  \sin \phi_1^4}  \le D.
\end{align*}
Because $\theta^4_1 < \pi/2$ and 
$\sin \theta_1^4  \sin \phi_1^4 < 1$, we have
\begin{align*}
\sin \theta_1^4  \sin \phi_1^4 \geq \frac{6}{D} =: M, \qquad
  0 < M < 1.
\end{align*}
By Lemma~\ref{lem59} with $j = 4$ and $p = 1$, setting 
$\gamma (M) := \pi - \sin^{-1} M$, we have 
$\frac{\pi}{2} < \gamma(M) < \pi$, and the MAC
with $\gamma(M)$ is satisfied on $F_2$, $F_3$, $F_4$, and 
$\psi^{2,4}$, $\psi^{3,4} \, \revise{\le} \, \gamma(M)$.

Note that $|T| = \frac{1}{6}\alpha_1\alpha_2\alpha_3
\sin \theta_1^3 \sin \phi_1^3$, and we have
\begin{align*}
  \frac{H_T}{h_T} = \frac{\alpha_1 \alpha_2 \alpha_3}{|T|}
   = \frac{6}{\sin \theta_1^3 \sin \phi_1^3}  \le D.
\end{align*}
Thus, by Lemma~\ref{lem59} with $j = 3$ and $p = 1$,
we find that $\psi^{2,3} \, \revise{\le} \, \gamma(M)$.

Because $|\overline{P_3 P_4}| < |\overline{P_1 P_4}| +
|\overline{P_1P_3}| \leq 2 \alpha_3$
on $F_2 = \triangle P_1 P_3 P_4$
and $|\overline{P_2P_3}| \leq \alpha_1$, we note that
\begin{align*}
  |T| = \frac{1}{6} \alpha_2  |\overline{P_2 P_3}| |\overline{P_3 P_4}|
    \sin \theta^1_3 \sin \phi_3^1 
   < \frac{1}{3} \alpha_1 \alpha_2 \alpha_3 \sin \theta^1_3 \sin \phi_3^1.
\end{align*}
Thus, we have
\begin{align*}
  D \ge \frac{H_K}{h_K} > \frac{3}{\sin \theta^1_3 \sin \phi_3^1}
  \quad \text{ and } \quad
  \sin \theta^1_3 \sin \phi_3^1 > \frac{3}{D} = \frac{M}{2}.
\end{align*}
From Lemma~\ref{lem59}, setting $\gamma (M/2) := \pi - \sin^{-1} (M/2)$,
we have $\frac{\pi}{2} < \gamma(M/2) < \pi$ and 
MAC with $\gamma(M/2)$ is satisfied on $F_1$,
and $\psi^{2,1}$, $\psi^{4,1} \, \revise{\le} \,  \gamma(M/2)$.

The final thing to prove is the MAC for $\psi^{1,3}$.
From the cosine rule \eqref{semi17d}, we have
\begin{align*}
\cos \psi^{1,3} = \sin \psi^{3,4} \sin \psi^{4,1}
   \cos \theta_2^{4} - \cos \psi^{3,4} \cos \psi^{4,1}.
\end{align*}
By the definition of Type~1, the angle $\theta_2^4$ is the minimum angle
of $F_4 = \triangle P_1 P_2 P_3$, and therefore, we have
\begin{align*}
   \cos \theta_2^4 \ge \frac{1}{2}, \;
\sin \psi^{{3},{4}} \sin \psi^{{4},{1}} \cos \theta^{4}_2 > 0,
\quad\text{ and }
 \cos \psi^{1,3} > - \cos \psi^{3,4} \cos \psi^{4,1}.
\end{align*}
From the above argument, we have $\sin \psi^{3,4} > M$,
$\sin \psi^{4,1} > M/2$, and
\begin{align*}
\cos \psi^{1,3} & > - \cos \psi^{3,4} \cos \psi^{4,1} 
 \ge - | \cos \psi^{3,4} | | \cos \psi^{4,1} | \\
 & = - \sqrt{1 - \sin^2 \psi^{3,4}} \sqrt{1 - \sin^2 \psi^{4,1}}
 > - \sqrt{1 - M^2} \sqrt{1 - \frac{M^2}{4}} > -1.
\end{align*}
Therefore, we conclude that
\begin{align*}
\psi^{1,3} < \cos^{-1}
  \left(- \sqrt{1 - M^2} \sqrt{1 - \frac{M^2}{4}}\right) < \pi,
\end{align*}
and $T$ satisfies the MAC with
\begin{align*}
  \gmax := \max\left\{\gamma(M/2),  \cos^{-1}
  \left(- \sqrt{1 - M^2} \sqrt{1 - \frac{M^2}{4}}\right)
   \right\}.
\end{align*}

\subsection{Type~2: Proof of ``MAC implies \eqref{equiv-cond}''}
First, we suppose that $T$ is of Type~2 and satisfies the MAC
with $\gmax \, \revise{\in [\pi/3, \pi)}$. 
  The proof is very similar to
that described in Section~\ref{Case1}.

By the definition of Type~2,
$\alpha_3 = |\overline{P_1P_4}| \le |\overline{P_2P_4}|$.
Because
\begin{align*}
    |T| = \frac{1}{6} \alpha_1 \alpha_2 \alpha_3
   \sin \theta_2^4  \sin \phi_1^4 =
    \frac{1}{6} \alpha_1 \alpha_2 |\overline{P_2P_4}|
    \sin \theta_2^4  \sin \phi_2^4,
\end{align*}
we have
 \begin{align}
   \frac{H_T}{h_T} = \frac{\alpha_1 \alpha_2 \alpha_3}{|T|}
     = \frac{6}{\sin \theta_2^4  \sin \phi_1^4}
  \le \frac{6}{\sin \theta_2^4  \sin \phi_2^4}.
  \label{Case3}
\end{align}
From the definition of Type~2, we realize that
$\theta_1^4 \le \theta_2^4 \le \theta_3^4$ on $F_4$, 
$\theta_2^3 \le \theta_1^3 \le \theta_4^3$ on $F_3$, and
$\theta_2^1$ is not the minimum angle of $F_1$.
Thus, it follows from Lemma \ref{lem53} that
\begin{align*}
 \frac{\pi - \gmax}{2} \leq \theta_2^4, \, \theta_1^3, 
  \, \theta_2^1\le \gmax,
  \quad  \sin \theta_2^4, \, \sin \theta_1^3, \,
  \sin \theta_2^1 \ge C_1.
\end{align*}
Additionally, we may apply Lemma~\ref{lem56} to $\theta_2^4$ and $F_4$,
and find that either $\psi^{1,4} \ge \delta$ or
$\psi^{3,4} \ge \delta$, where $\delta = \delta(\gmax)$ is defined by
\eqref{def-delta}.

Suppose that $\psi^{3,4} \ge \delta$.  Using the same argument as in
Section~\ref{Case1}, we have
\begin{align*}
  \sin \phi_1^4 = \sin \theta_1^3 \sin \psi^{3,4}
  \ge C_0 \sin \theta_1^3 \ge C_0C_1.
\end{align*}

Next, suppose that $\psi^{1,4} \ge \delta$.  We have 
\begin{align*}
  \sin \phi_2^4 = \sin \theta_2^1 \sin \psi^{1,4}
  \ge C_0 \sin \theta_2^1 \ge C_0C_1.
\end{align*}

Combining these results with \eqref{Case3}, we obtain
\begin{align*}
  \frac{H_T}{h_T}  \le \frac{6}{C_0 C_1^2} =: D,
\end{align*}
that is, \eqref{equiv-cond} holds.

\subsection{Type~2: Proof of ``\eqref{equiv-cond} implies MAC''}
Finally, we suppose that $T$ is of Type~2 and 
\begin{align*}
  \frac{H_T}{h_T} = \frac{\alpha_1 \alpha_2 \alpha_3}{|T|}
   = \frac{6}{\sin \theta_2^4  \sin \phi_1^4} \le D, \quad
   \sin \theta_2^4  \sin \phi_1^4 \ge \frac{6}{D} =: M. 
\end{align*}
The proof is very similar to that described in Section~\ref{Case2}.
By Lemma~\ref{lem59} with $j = 4$ and $p = 1$, setting 
$\gamma (M) := \pi - \sin^{-1} M$, 
the MAC with $\gamma(M)$ is satisfied on
$F_2$, $F_3$, $F_4$, and 
$\psi^{2,4}$, $\psi^{3,4} \, \revise{\le} \, \gamma(M)$.

Because $|\overline{P_2 P_4}| \le \alpha_1$, we have
\begin{align*}
  |T| = \frac{1}{6} |\overline{P_2 P_3}|
  |\overline{P_2 P_4}| |\overline{P_1 P_4}|
    \sin \theta_2^1 \sin \phi_4^1 
   \le \frac{1}{6} \alpha_1 \alpha_2 \alpha_3 \sin \theta_2^1
    \sin \phi_4^1.
\end{align*}
This yields
\begin{align*}
  D \ge \frac{H_T}{h_T} \ge \frac{6}{\sin \theta_2^1 \sin \phi_4^1}
  \quad \text{ and } \quad
  \sin \theta_2^1 \sin \phi_4^1 \ge \frac{6}{D} = M,
\end{align*}
and, by Lemma~\ref{lem59} with $j = 1$ and $p = 4$, we find 
that the MAC with $\gamma(M)$ is satisfied on $F_1$,
and $\psi^{1,2}$, $\psi^{1,3} \, \revise{\le} \, \gamma(M)$.

The final thing to prove is the MAC for
$\psi^{1,4}$ and $\psi^{2,3}$.  
By the cosine rule \eqref{semi17d} with $j=2$, we have
\begin{align*}
  \cos \psi^{1,4} & = \sin \psi^{1,3} \sin \psi^{4,3}
   \cos \theta_2^3 - \cos \psi^{1,3} \cos \psi^{4,3}, \\
   \cos \psi^{2,3} & = \sin \psi^{2,4} \sin \psi^{3,4}
   \cos \theta_1^4 - \cos \psi^{2,4} \cos \psi^{3,4}.
\end{align*}
By the definition of Type 2,
$\theta_2^3$ and $\theta_1^4$ are the minimum angles of $F_3$
and $F_4$, respectively.  Therefore, we have
$\cos \theta_2^3$, $\cos \theta_1^4 \ge \frac{1}{3}$ and thus
\begin{align*}
 \cos \psi^{1,4} > - \cos \psi^{1,3} \cos \psi^{3,4}, \qquad
 \cos \psi^{2,3} > - \cos \psi^{2,4} \cos \psi^{3,4}.
\end{align*}
Because $\sin \psi^{1,3}$, $\sin \psi^{2,4}$, 
$\sin \psi^{3,4} > M$, we find that
\begin{align*}
   \cos \psi^{1,4} & > - \cos \psi^{1,3} \cos \psi^{3,4} 
    \ge - \sqrt{1 - \sin^2 \psi^{1,3}} \sqrt{1 - \sin^2 \psi^{3,4}}
    > M^2 - 1, \\
   \cos \psi^{2,3} & > M^2 - 1.
\end{align*}
Therefore, we conclude that
$\psi^{1,4}$, $\psi^{2,3} < \cos^{-1}(M^2 - 1) < \pi$,
and $T$ satisfies the MAC with
\begin{align*}
  \gmax := \max\left\{\gamma(M), \cos^{-1}(M^2 - 1)
   \right\}.
\end{align*}

\section{Concluding remarks}
The equivalence between the maximum angle condition and the boundedness
of $R_T/h_T$ (and $H_T/h_T$) \revise{has been} established.  Because the
ratio $R_T/h_T$ appears in many error estimations, the equivalence
relationship is very valuable in the mathematical theory of finite
element methods. Additionally, because $h_T$ and $R_T$ can be easily
computed, \revise{this will} hopefully \revise{be} very useful in many
practical computations such as mesh generation and adaptive mesh refinement.

Mathematically, an interesting and challenging problem is to
extend this result to the case of $d$-simplices $(d \ge 4)$.
To do this, we need to develop the theory of $d$-simplex geometry 
further to obtain deeper insight\revise{s}.

\vspace{0.5cm}
\noindent
\textbf{Acknowledgments}
The second and fourth authors were supported by JSPS KAKENHI Grant
Number 20H01820.
\revise{The authors thank the anonymous referees for their
valuable comments.}

\end{document}